\documentclass[a4paper,12pt]{amsart}
\usepackage{amsfonts}
\usepackage{amssymb}
\usepackage{ifthen}
\usepackage{amscd}
\usepackage{amsxtra}
\usepackage{graphicx}
\usepackage{color}
\nonstopmode \numberwithin{equation}{section}
\newtheorem{thm}{Theorem}
\newtheorem{lem}{Lemma}
\newtheorem{cor}{Corollary}[section]

\newtheorem{cl}{Claim}
\newtheorem{ca}{Case}
\newtheorem{sca}{Subcase}
\newtheorem{scl}{Subclaim}
\newtheorem{conj}{Conjecture}

\theoremstyle{definition}
\newtheorem{defn}{Definition}

\newtheorem{op}[equation]{Open Problem}
\newtheorem{ques}[equation]{Question}
\newtheorem{rem}{Remark}[section]
\newtheorem{exam}[equation]{Example}

\newcounter {own}
\def\theown {\thesection       .\arabic{own}}

\newenvironment{pf}[1][]{%
 \vskip 3mm
 \noindent
 \ifthenelse{\equal{#1}{}}%
  {{\slshape Proof. }}%
  {{\slshape #1.} }%
 }%
{\qed\bigskip}

\newcounter{alphabet}
\newcounter{tmp}
\newenvironment{Thm}[1][]{\refstepcounter{alphabet}%
\bigskip%
\noindent%
{\bf Theorem \Alph{alphabet}}%
\ifthenelse{\equal{#1}{}}{}{ (#1)}%
{\bf .} \itshape}{\vskip 8pt}

\makeatletter
\newcommand{\Ref}[1]{\@ifundefined{r@#1}{}{\setcounter{tmp}{\ref{#1}}\Alph{tmp}}}
\makeatother

\newenvironment{Lem}[1][]{\refstepcounter{alphabet}%
\bigskip%
\noindent%
{\bf Lemma \Alph{alphabet}}%
{\bf .} \itshape}{\vskip 8pt}

\newcommand{\ID}{{\mathbb D}}

\def\be{\begin{equation}}
\def\ee{\end{equation}}

\newcommand{\bee}{\begin{enumerate}}
\newcommand{\eee}{\end{enumerate}}

\newcommand{\blem}{\begin{lem}}
\newcommand{\elem}{\end{lem}}
\newcommand{\bthm}{\begin{thm}}
\newcommand{\ethm}{\end{thm}}
\newcommand{\bcor}{\begin{cor}}
\newcommand{\ecor}{\end{cor}}
\newcommand{\beg}{\begin{exam}}
\newcommand{\eeg}{\end{exam}}
\newcommand{\begs}{\begin{examples}}
\newcommand{\eegs}{\end{examples}}
\newcommand{\bdefe}{\begin{defn}}
\newcommand{\edefe}{\end{defn}}
\newcommand{\bprob}{\begin{prob}}
\newcommand{\eprob}{\end{prob}}
\newcommand{\bques}{\begin{ques}}
\newcommand{\eques}{\end{ques}}
\newcommand{\bei}{\begin{itemize}}
\newcommand{\eei}{\end{itemize}}
\newcommand{\bcon}{\begin{conj}}
\newcommand{\econ}{\end{conj}}
\newcommand{\bop}{\begin{op}}
\newcommand{\eop}{\end{op}}

\newcommand{\bca}{\begin{ca}}
\newcommand{\eca}{\end{ca}}
\newcommand{\bsca}{\begin{sca}}
\newcommand{\esca}{\end{sca}}

\newcommand{\bcl}{\begin{cl}}
\newcommand{\ecl}{\end{cl}}

\newcommand{\bscl}{\begin{scl}}
\newcommand{\escl}{\end{scl}}

\newcommand{\bcons}{\begin{conjs}}
\newcommand{\econs}{\end{conjs}}
\newcommand{\bprop}{\begin{propo}}
\newcommand{\eprop}{\end{propo}}
\newcommand{\br}{\begin{rem}}
\newcommand{\er}{\end{rem}}
\newcommand{\brs}{\begin{rems}}
\newcommand{\ers}{\end{rems}}
\newcommand{\bo}{\begin{obser}}
\newcommand{\eo}{\end{obser}}
\newcommand{\bos}{\begin{obsers}}
\newcommand{\eos}{\end{obsers}}
\newcommand{\bpf}{\begin{pf}}
\newcommand{\epf}{\end{pf}}
\newcommand{\ba}{\begin{array}}
\newcommand{\ea}{\end{array}}
\newcommand{\beq}{\begin{eqnarray}}
\newcommand{\beqq}{\begin{eqnarray*}}
\newcommand{\eeq}{\end{eqnarray}}
\newcommand{\eeqq}{\end{eqnarray*}}

\newcommand{\ds}{\displaystyle}

\newcounter{minutes}
\divide\time by 60
\newcounter{hours}
\multiply\time by 60 \addtocounter{minutes}{-\time}

\begin{document}

\bibliographystyle{amsplain}
\title []
{Some properties of  a class of elliptic partial differential
operators}

\def\thefootnote{}
\footnotetext{ \texttt{\tiny File:~\jobname .tex,
          printed: \number\day-\number\month-\number\year,
          \thehours.\ifnum\theminutes<10{0}\fi\theminutes}
} \makeatletter\def\thefootnote{\@arabic\c@footnote}\makeatother

\author{Shaolin Chen}
\address{Shaolin Chen, Department of Mathematics and Computational
Science, Hengyang Normal University, Hengyang, Hunan 421008,
People's Republic of China.} \email{mathechen@126.com}

\author{Matti Vuorinen}
\address{Matti Vuorinen,
Department of Mathematics and Statistics,
 University of Turku, Turku 20014,
Finland. }\email{vuorinen@utu.fi}





\begin{abstract} We prove Schwarz-Pick type estimates and
coefficient estimates for a class of functions induced by the
elliptic partial differential operators. Then we apply these results
to obtain a Landau type theorem.
\end{abstract}


\maketitle \pagestyle{myheadings} \markboth{ Sh. Chen and M.
Vuorinen}{Differential operators}

\section{Introduction and main results }\label{csw-sec1}
Let  $\mathbb{C}$ be the complex plane. For $a\in\mathbb{C}$, let
$r>0$ and $\ID(a,r)=\{z:\, |z-a|<r\}$. In particular, we use
$\mathbb{D}_r$ to denote the disk $\mathbb{D}(0,r)$ and
$\mathbb{D}$, the open unit
disk $\ID_1$. 

For a real $2\times2$ matrix,
 we will consider the matrix norm $\|A\|=\sup\{|Az|:~|z|=1\}$ and
the matrix function $l(A)=\inf\{|Az|:~|z|=1\}$. For
$z=x+iy\in\mathbb{C}$ with $x$ and $y$ real, we denote the {\it
complex differential operators}
$$\frac{\partial}{\partial
z}=\frac{1}{2}\left(\frac{\partial}{\partial
x}-i\frac{\partial}{\partial
y}\right)~\mbox{and}~\frac{\partial}{\partial
\overline{z}}=\frac{1}{2}\left(\frac{\partial}{\partial
x}+i\frac{\partial}{\partial y}\right).
$$ If we denote the formal derivative of $f=u+iv$ by

$$D_{f}=\left(\begin{array}{cccc}
\ds u_{x}\;~~ u_{y}\\[2mm]
\ds v_{x}\;~~ v_{y}
\end{array}\right),
$$ then $\|D_{f}\|=|f_{z}|+|f_{\overline{z}}|$ and
$l(D_{f})=\big||f_{z}|-|f_{\overline{z}}|\big|,$ where $u$, $v$ are
real functions, $f_{z}=\partial f/\partial z$ and
$f_{\overline{z}}=\partial f/\partial \overline{z}$. Throughout this
paper, we denote by $ \mathcal{C}^{n}(\mathbb{D})$ the set of all
$n$-times continuously differentiable complex-valued functions in
$\mathbb{D}$, where $n\in\{1,2,\cdots\}$.

For $\alpha\in\mathbb{R}$ and $z\in\mathbb{D}$, let

$$T_{\alpha}=-\frac{\alpha^{2}}{4}(1-|z|^{2})^{-\alpha-1}+\frac{\alpha}{2}(1-|z|^{2})^{-\alpha-1}
\left(z\frac{\partial}{\partial
z}+\overline{z}\frac{\partial}{\partial
\overline{z}}\right)+\frac{1}{4}(1-|z|^{2})^{-\alpha}\Delta$$ be the
{\it second order elliptic partial differential operator}, where
$\Delta$ is the usual complex {\it Laplacian operator}
$$\Delta:=4\frac{\partial^{2}}{\partial z\partial
\overline{z}}=\frac{\partial^{2}}{\partial
x^{^{2}}}+\frac{\partial^{2}}{\partial y^{^{2}}}.
$$

We consider the {\it Dirichlet boundary value problem} of
distributional sense as follows

\be\label{eq-1}
\begin{cases}
\displaystyle T_{\alpha}(f)=0
& \mbox{in } \mathbb{D},\\
\displaystyle~ f=f^{\ast} &\mbox{on } \partial\mathbb{D}.
\end{cases}
\ee Here, the boundary data
$f^{\ast}\in\mathfrak{D}'(\partial\mathbb{D})$ is a {\it
distribution} on the boundary  $\partial\mathbb{D}$ of $\mathbb{D}$,
and the boundary condition in (\ref{eq-1}) is interpreted in the
distributional sense that $f_{r}\rightarrow f^{\ast}$ in
$\mathfrak{D}'(\partial\mathbb{D})$ as $r\rightarrow1-$, where
\be\label{eq-0.1}f_{r}(e^{i\theta})=f(re^{i\theta}),~e^{i\theta}\in\partial\mathbb{D},\ee
for $r\in[0,1)$ (see \cite{O}).

In \cite{O}, Olofsson proved  that, for parameter values
$\alpha>-1$, a function $f\in\mathcal{C}^{2}(\mathbb{D})$ satisfies
(\ref{eq-1}) if and only if it has the form of a {\it Poisson type
integral}
\be\label{eq-0.2}f(z)=\frac{1}{2\pi}\int_{0}^{2\pi}K_{\alpha}(ze^{-i\tau})f^{\ast}(e^{i\tau})d\tau,~\mbox{for}~z\in\mathbb{D},\ee
where
$$K_{\alpha}(z)=c_{\alpha}\frac{(1-|z|^{2})^{\alpha+1}}{|1-z|^{\alpha+2}},$$
 $c_{\alpha}=\big(\Gamma(\alpha/2+1)\big)^{2}/\Gamma(1+\alpha)$ and
 $\Gamma(s)=\int_{0}^{\infty}t^{s-1}e^{-t}dt$ for $s>0$ is the
 standard Gamma function. If we take $\alpha=2(n-1)$, then $f$ is
 {\it polyharmonic} (or {\it $n$-harmonic}), where $n\in\{1,2,\ldots\}$
 (cf. \cite{AA,AP,AH,E,O1,O2,O3,P1}). Furthermore, Borichev and Hedenmalm \cite{AH} proved that
 $$(1-|z|^{2})^{n}\Delta^{n}=4(1-|z|^{2})T_{0}\circ4(1-|z|^{2})^{2}T_{2}\circ\cdots\circ4(1-|z|^{2})^{n}T_{2(n-1)}.$$
 In particular, if $\alpha=0$, then $f$ is
 harmonic.

 For $a, b, c\in\mathbb{R}$ with $c\neq0, -1, -2, \ldots,$ the {\it
hypergeometric} function is defined by the power series
$$F(a,b;c;x)=\sum_{n=0}^{\infty}\frac{(a)_{n}(b)_{n}}{(c)_{n}}\frac{x^{n}}{n!},~|x|<1,$$
where $(a)_{0}=1$ and $(a)_{n}=a(a+1)\cdots(a+n-1)$ for $n=1, 2,
\ldots$ are the {\it Pochhammer} symbols. Obviously, for  $n=0, 1,
2, \ldots$, $(a)_{n}=\Gamma(a+n)/\Gamma(a).$ In particular, for $a,
b, c>0$ and $a+b<c$, we have (cf. \cite{ABRVV,QVV}) \be\label{eq-y}
F(a,b;c;1)=\lim_{x\rightarrow1}
F(a,b;c;x)=\frac{\Gamma(c)\Gamma(c-a-b)}{\Gamma(c-a)\Gamma(c-b)}<\infty.\ee

 For
$\alpha=0$, Heinz \cite{HZ} and Colonna \cite{Co} proved the
following Schwarz-Pick type estimates on planar harmonic functions,
which are the following. For the related discussions on this topic,
 see \cite{H1,CPW2,CR,KV,P0}.


\begin{Thm}{\rm (\cite[Lemma]{HZ})}\label{HZ}
Let $f$ be a harmonic function of $\mathbb{D}$ into $\mathbb{D}$
with $f(0)=0.$ Then, for $z\in\mathbb{D}$,
$$|f(z)|\leq\frac{4}{\pi}\arctan|z|.
$$
This estimate is sharp.
\end{Thm}

\begin{Thm}{\rm (\cite[Theorems 3 \mbox{and} 4]{Co})}\label{Co}
Let $f$ be a harmonic function of $\mathbb{D}$ into $\mathbb{D}$.
Then, for $z\in\mathbb{D}$,
$$\|D_{f}(z)\|\leq\frac{4}{\pi}\frac{1}{1-|z|^{2}}.
$$
This estimate is sharp, and all the extremal functions are
$$f(z)=\frac{2\gamma }{\pi}\arg \left (
\frac{1+\psi(z)}{1-\psi(z)}\right), $$ where $|\gamma|=1$ and $\psi$
is a conformal automorphism of $\mathbb{D}$.
\end{Thm}

For $\alpha>-1$, we establish the following Schwarz-Pick type
estimate on the solutions to (\ref{eq-1}).

\begin{thm}\label{thm-2}
For $\alpha>-1$, let $f\in\mathcal{C}^{2}(\mathbb{D})$ satisfy {\rm
(\ref{eq-1})} and $\sup_{z\in\mathbb{D}}|f(z)|\leq M$, where $M$ is
a positive constant. Then, for $z\in\mathbb{D},$
\be\label{eq-0.1x}\left|f(z)-\frac{(1-|z|)^{\alpha+1}}{1+|z|}f(0)\right|\leq
M\left[\frac{1}{2\pi}\int_{0}^{2\pi}K_{\alpha}(ze^{-it})dt-\frac{(1-|z|)^{\alpha+1}}{1+|z|}K_{\alpha}(0)\right]\ee
and

\be\label{eq-0.12}\|D_{f}(z)\|\leq\frac{M\mathcal{M}_{\alpha}(|z|)
\big[2+\alpha+(4+3\alpha)|z|\big]}{1-|z|^{2}}\leq\frac{M
\big[2+\alpha+(4+3\alpha)|z|\big]}{1-|z|^{2}},\ee where
\be\label{eq-11}
\mathcal{M}_{\alpha}(r)=\frac{1}{2\pi}\int_{0}^{2\pi}K_{\alpha}(re^{i\epsilon})d\epsilon=
\frac{\left[\Gamma\big(1+\frac{\alpha}{2}\big)\right]^{2}}{\Gamma(1+\alpha)}
F\left(-\frac{\alpha}{2},-\frac{\alpha}{2};1;r^{2}\right),~r\in[0,1).\ee
\end{thm}

Let $f$ be a  harmonic mapping of $\overline{\mathbb{D}}$ onto
$\overline{\mathbb{D}}$ with $f(0)=0$. In \cite{HZ}, Heinz showed
that, for $\theta\in[0,2\pi]$,
$$\|D_{f}(e^{i\theta})\|\geq\frac{2}{\pi}.$$ For the extensive discussion
on the Heinz's inequality for real harmonic functions in high
dimension, see \cite{K}.

By using Theorem \ref{thm-2}, we  get a Heinz type inequality on
$\partial\mathbb{D}$ as follows.

\begin{thm}\label{thm-3}
For $\alpha\geq0$, let $f\in\mathcal{C}^{2}(\overline{\mathbb{D}})$
satisfying {\rm (\ref{eq-1})}. Suppose that $f(0)=0$,
$f(\overline{\mathbb{D}})=\overline{\mathbb{D}}$ and
$f(\partial\mathbb{D})=\partial\mathbb{D}.$

\begin{enumerate}
\item[\rm(a)] If $\alpha=0$, then, for $\theta\in[0,2\pi]$,
$$\|D_{f}(e^{i\theta})\|\geq\frac{2}{\pi};$$

\item[{\rm(b)}] If $\alpha>0$, then, for $\theta\in[0,2\pi]$,
$$\|D_{f}(e^{i\theta})\|\geq \lim_{r\rightarrow1-}\frac{d}{dr}\mathcal{M}_{\alpha}(r)=\frac{\alpha}{2},$$
where $\mathcal{M}_{\alpha}(r)$ is given by {\rm (\ref{eq-11})}.
\end{enumerate}
\end{thm}

The following result is the homogeneous expansion of solutions to
{\rm (\ref{eq-1})}.

\begin{Thm} $($\cite[Theorem 2.2]{O}$)$\label{ThmA}
Let $\alpha\in\mathbb{R}$ and $f\in\mathcal{C}^{2}(\mathbb{D})$.
Then $f$ satisfies {\rm (\ref{eq-1})} if and only if it has a series
expansion of the form \beq\label{eq-2}
f(z)&=&\sum_{k=0}^{\infty}c_{k}F\left(-\frac{\alpha}{2},k-\frac{\alpha}{2};k+1;|z|^{2}\right)z^{k}\\
\nonumber &&+
\sum_{k=1}^{\infty}c_{-k}F\left(-\frac{\alpha}{2},k-\frac{\alpha}{2};k+1;|z|^{2}\right)\overline{z}^{k},~z\in\mathbb{D},
\eeq for some sequence $\{c_{k}\}_{k=-\infty}^{\infty}$ of complex
numbers satisfying
\be\label{eq-3}\lim_{|k|\rightarrow\infty}\sup|c_{k}|^{\frac{1}{|k|}}\leq1.\ee
In particular, the expansion {\rm (\ref{eq-2})}, subject to {\rm
(\ref{eq-3})}, converges in $\mathcal{C}^{\infty}(\mathbb{D})$, and
every solution $f$ of {\rm (\ref{eq-1})} is
$\mathcal{C}^{\infty}$-smooth in $\mathbb{D}.$
\end{Thm}

For $\alpha=0$, there are numerous discussions on coefficient
estimates of harmonic mappings in the literature,
 see for example \cite{CPW1,CPW2,CPW3,CPW4, Du, LZ, SW}. We
 investigate the problem of coefficient estimates  on the solutions to
 (\ref{eq-1}) as follows.

\begin{thm}\label{thm-1}
For $\alpha>-1$, let $f\in\mathcal{C}^{2}(\mathbb{D})$ be a solution
to {\rm (\ref{eq-1})} with the series expansion of the form {\rm
(\ref{eq-2})} and $\sup_{z\in\mathbb{D}}|f(z)|\leq M$, where $M$ is
a positive constant. Then, for $k\in\{1,2,\ldots\}$,
\be\label{eq-4}\left|c_{k}F\Big(-\frac{\alpha}{2},k-\frac{\alpha}{2};k+1;1\Big)\right|+
\left|c_{-k}F\Big(-\frac{\alpha}{2},k-\frac{\alpha}{2};k+1;1\Big)\right|\leq\frac{4M}{\pi}\ee
and
$$\Big|c_{0}F\left(-\frac{\alpha}{2},-\frac{\alpha}{2};1;1\right)\Big|\leq
M.$$ In particular, if $\alpha=0$, then the estimate of {\rm
(\ref{eq-4})} is sharp and all the extreme functions are
$$f_{k}(z)=\frac{2\varepsilon M}{\pi}{\rm
Im}\left(\log\frac{1+\vartheta z^{k}}{1-\vartheta z^{k}}\right),$$
where $|\varepsilon|=|\vartheta |=1$.
\end{thm}

The following result easily follows from Theorem \ref{thm-1} and
\cite[Proposition 1.4]{O}.

\begin{cor}
For $\alpha>-1$, let $f\in\mathcal{C}^{2}(\mathbb{D})$ be a solution
to {\rm (\ref{eq-1})} with the series expansion of the form {\rm
(\ref{eq-2})} and $\sup_{z\in\mathbb{D}}|f(z)|\leq M$, where $M$ is
a positive constant. Then, for $k\in\{1,2,\ldots\}$,
$$|c_{k}|+|c_{-k}|\leq\frac{4M\Gamma\left(1+\frac{\alpha}{2}\right)\Gamma
\left(k+1+\frac{\alpha}{2}\right)}{k!\Gamma(\alpha+1)\pi}.$$

\end{cor}

For $p\in(0,\infty]$, the {\it Hardy space $\mathcal{H}^{p}$}
 consists of those functions $f:\ \mathbb{D}\rightarrow\mathbb{C}$
 such that $f$ is measurable, $M_{p}(r,f)$ exists for all $r\in(0,1)$ and  $ \|f\|_{p}<\infty$, where
$$
\|f\|_{p}=
\begin{cases}
\displaystyle\sup_{0<r<1}M_{p}(r,f),
& \mbox{ if } p\in(0,\infty),\\
\displaystyle\sup_{z\in\mathbb{D}}|f(z)|, &\mbox{ if }\, p=\infty,
\end{cases}$$ 
and
$$M_{p}^{p}(r,f)=\frac{1}{2\pi}\int_{0}^{2\pi}|f(re^{i\theta})|^{p}\,d\theta.
$$ It is not difficult to know that all bounded measurable functions
belong to $\mathcal{H}^{p}$.

The classical theorem of Landau  shows that there is a $\rho
=\frac{1}{M+\sqrt{M^{2}-1}}$ such that every function $f$, analytic
in $\mathbb{D}$ with $f(0)=f'(0)-1=0$ and $|f(z)|<M$ in
$\mathbb{D}$, is univalent in the disk $\mathbb{D}_\rho $ and in
addition, the range $f(\mathbb{D}_\rho )$ contains a disk of radius
$M\rho ^2$ (see \cite{L}), where $M\geq1$ is a constant. Recently,
many authors considered Landau type theorem for planar harmonic
mappings (see \cite{HG,CPW1,CPW2,CPW3}),
 biharmonic mappings (see \cite{AA}) and polyharmonic mappings (see \cite{CPW4}).

The last main aim of this paper is to establish a landau type
theorem for a more general class of functions. Moreover, we do not
need the condition of bounded functions as in the classical Landau
Theorem. Applying Theorems \ref{thm-2} and \ref{thm-1}, we get the
following Landau type theorem for a class of functions $f\in
\mathcal{H}^{p}$ satisfying {\rm (\ref{eq-1})}.

\begin{thm}\label{thm-4}
For $\alpha\in(-1,0]$, let $f\in\mathcal{C}^{2}(\mathbb{D})$ be a
solution to  {\rm (\ref{eq-1})} satisfying
$f(0)=|J_{f}(0)|-\lambda=0$ and $f\in\mathcal{H}^{p}$, where
$\lambda$ is a positive constant and $J_{f}$ is the Jacobian of $f$.
Then $f$ is univalent in $\mathbb{D}_{\gamma_{0}\rho_{0}}$, where
$\rho_{0}$ satisfies the following equation
$$\frac{\lambda}{M^{\ast}(2+\alpha)}-\frac{4M^{\ast}\rho_{0}}{\pi}\Bigg[\frac{2-\rho_{0}}{(1-\rho_{0})^{2}}
+\frac{2\rho_{0}}{(1-\rho_{0})(1-\rho_{0}^{2})^{2}}\Bigg]=0,$$ where
$\mu(\gamma)=\big(1+\gamma\big)^{\frac{\alpha+1}{p}}/
\left[\gamma\big(1-\gamma\big)^{\frac{1}{p}}\right]$,
$\mu(\gamma_{0})=\min_{0<\gamma<1}\mu(\gamma)$ and
$M^{\ast}=c_{\alpha}^{\frac{1}{p}}\|f\|_{p}\mu(\gamma_{0}).$
Moreover, $f\big(\mathbb{D}_{\gamma_{0}\rho_{0}}\big)$ contains a
univalent disk $\mathbb{D}_{\gamma_{0}R_{0}}$ with
$$R_{0}\geq\frac{2\rho_{0}}{3}\left[\frac{\lambda}{M^{\ast}(2+\alpha)}-\frac{M^{\ast}\rho_{0}(2-\rho_{0})}{\pi(1-\rho_{0})^{2}}\right].
$$
\end{thm}

We remark that Theorem \ref{thm-4} is a generalization of
\cite[Theorem 2]{HG} and \cite[Theorem 5]{CPW5}.

The proofs of Theorems \ref{thm-2}, \ref{thm-3} and  \ref{thm-1}
 will be presented in Section \ref{csw-sec2}, and the
proof of Theorem \ref{thm-4} will be given in Section
\ref{csw-sec3}.

\section{ Schwarz-Pick type estimates and  coefficient estimates }\label{csw-sec2}

\subsection*{Proof of Theorem \ref{thm-2}} We first prove (\ref{eq-0.1x}).  By the assumption, we
see that $f_{r}\rightarrow f^{\ast}$ in
$\mathfrak{D}'(\partial\mathbb{D})$ as $r\rightarrow1-$, where
$f_{r}$ is given by (\ref{eq-0.1}) for $r\in[0,1).$
 By (\ref{eq-0.2}), for $z=re^{i\theta}\in\mathbb{D}$,
we have

\begin{eqnarray*}
&&\left|f(z)-\frac{(1-|z|)^{\alpha+1}}{1+|z|}f(0)\right|\\
&=&\frac{1}{2\pi}\left|\int_{0}^{2\pi}K_{\alpha}(ze^{-it})f^{\ast}(e^{it})dt-\frac{(1-|z|)^{\alpha+1}}{1+|z|}
\int_{0}^{2\pi}K_{\alpha}(0)f^{\ast}(e^{it})dt \right|\\
&=&\frac{1}{2\pi}\left|\int_{0}^{2\pi}\Big(K_{\alpha}(ze^{-it})-\frac{(1-|z|)^{\alpha+1}}{1+|z|}K_{\alpha}(0)\Big)f^{\ast}(e^{it})dt\right|\\
&\leq&\frac{1}{2\pi}\int_{0}^{2\pi}\Big(K_{\alpha}(ze^{-it})-\frac{(1-|z|)^{\alpha+1}}{1+|z|}K_{\alpha}(0)\Big)|f^{\ast}(e^{it})|dt\\
&\leq&M\left[\frac{1}{2\pi}\int_{0}^{2\pi}K_{\alpha}(ze^{-it})dt-\frac{(1-|z|)^{\alpha+1}}{1+|z|}K_{\alpha}(0)\right].
\end{eqnarray*}

Next we prove (\ref{eq-0.12}).  
By the proof of  \cite[Theorem 3.1]{O}, we observe that
\be\label{eq-0.11}\mathcal{M}_{\alpha}(r)=\frac{1}{2\pi}\int_{0}^{2\pi}K_{\alpha}(re^{i\epsilon})d\epsilon=
\frac{\left[\Gamma\big(1+\frac{\alpha}{2}\big)\right]^{2}}{\Gamma(1+\alpha)}
F\left(-\frac{\alpha}{2},-\frac{\alpha}{2};1;r^{2}\right)\ee and
$\mathcal{M}_{\alpha}(r)$ is increasing on $r\in[0,1)$ with
$$\lim_{r\rightarrow1-}\mathcal{M}_{\alpha}(r)=1.$$ By elementary calculations, for $z\in\mathbb{D}$, we have

$$\frac{\partial}{\partial
z}K_{\alpha}(ze^{-it})=c_{\alpha}\frac{(1-|z|^{2})^{\alpha}\left[\big(1+\frac{\alpha}{2}\big)e^{-it}(1-\overline{z}e^{it})(1-|z|^{2})-
(\alpha+1)\overline{z}|1-ze^{-it}|^{2}\right]}{|1-ze^{-it}|^{4+\alpha}}$$
and
$$\frac{\partial}{\partial
\overline{z}}K_{\alpha}(ze^{-it})=c_{\alpha}\frac{(1-|z|^{2})^{\alpha}\left[\big(1+\frac{\alpha}{2}\big)e^{it}(1-ze^{-it})(1-|z|^{2})-
(\alpha+1)z|1-ze^{-it}|^{2}\right]}{|1-ze^{-it}|^{4+\alpha}},$$
which, together with (\ref{eq-0.2}) and (\ref{eq-0.11}), imply that
\begin{eqnarray*}
\|D_{f}(z)\|&=&\left|\frac{1}{2\pi}\int_{0}^{2\pi}\frac{\partial}{\partial
z}K_{\alpha}(ze^{-it})f^{\ast}(e^{it})dt\right|+\left|\frac{1}{2\pi}\int_{0}^{2\pi}\frac{\partial}{\partial
\overline{z}}K_{\alpha}(ze^{-it})f^{\ast}(e^{it})dt\right|\\
&\leq&\frac{Mc_{\alpha}}{\pi}\int_{0}^{2\pi}\frac{(1-|z|^{2})^{\alpha}\big[(1+\alpha)|z||1-ze^{-it}|^{2}+
\big(1+\frac{\alpha}{2}\big)|1-ze^{-it}|(1-|z|^{2})\big]}{|1-ze^{-it}|^{4+\alpha}}dt\\
\end{eqnarray*}
\begin{eqnarray*}
&=&\frac{Mc_{\alpha}}{\pi} \left[\int_{0}^{2\pi}\frac{(1+\alpha)|z|(1-|z|^{2})^{\alpha}}{|1-ze^{-it}|^{2+\alpha}}dt+
\int_{0}^{2\pi}\frac{\big(1+\frac{\alpha}{2}\big)(1-|z|^{2})^{\alpha+1}}{|1-ze^{-it}|^{3+\alpha}}dt \right]\\
&\leq&\frac{Mc_{\alpha}}{\pi}\left[\frac{(1+\alpha)|z|}{(1-|z|^{2})}\int_{0}^{2\pi}\frac{(1-|z|^{2})^{\alpha+1}}{|1-ze^{-it}|^{2+\alpha}}dt
 +\frac{\big(1+\frac{\alpha}{2}\big)}{(1-|z|)}\int_{0}^{2\pi}\frac{(1-|z|^{2})^{\alpha+1}}{|1-ze^{-it}|^{2+\alpha}}dt\right] \\
&=&\left[\frac{2M(1+\alpha)|z|}{1-|z|^{2}}+\frac{M(2+\alpha)}{1-|z|}\right]\frac{1}{2\pi}\int_{0}^{2\pi}K_{\alpha}(ze^{-it})dt\\
&=&\left[\frac{2M(1+\alpha)|z|}{1-|z|^{2}}+\frac{M(2+\alpha)}{1-|z|}\right]\mathcal{M}_{\alpha}(|z|)\\
&=&
\frac{M\mathcal{M}_{\alpha}(|z|)\big[2+\alpha+(4+3\alpha)|z|\big]}{1-|z|^{2}}\\
&\leq&\frac{M\big[2+\alpha+(4+3\alpha)|z|\big]}{1-|z|^{2}}.
\end{eqnarray*}
The proof of this theorem is complete. \hfill $\Box$

\subsection*{Proof of Theorem \ref{thm-3}} Since (a) easily follows from the inequality (15) in \cite{HZ},
we only need to prove (b). Let $\alpha>0$. By Theorem \ref{thm-2}
(\ref{eq-0.1x}), we have

 \beq\label{eq-010}
\nonumber\frac{|f(e^{i\theta})-f(re^{i\theta})|}{1-r}&\geq&\frac{1-|f(re^{i\theta})|}{1-r}\\
&\geq&\frac{1-\frac{1}{2\pi}\int_{0}^{2\pi}K_{\alpha}(ze^{-it})dt+\frac{(1-|z|)^{\alpha+1}}{1+|z|}K_{\alpha}(0)}{1-r}.
\eeq where  $z=re^{i\theta}\in\mathbb{D}$ and $\theta\in[0,2\pi)$.

Applying \cite[Theorem 3.1]{O}, we get

$$\frac{1}{2\pi}\lim_{|z|\rightarrow1-}\int_{0}^{2\pi}K_{\alpha}(ze^{-it})dt=\lim_{r\rightarrow1-}\mathcal{M}_{\alpha}(r)=1,$$
which, together with L'Hopital's rule and (\ref{eq-010}), yield that

\begin{eqnarray*}
\|D_{f}(e^{i\theta})\|&\geq&\Bigg(\left|\frac{\partial
f(re^{i\theta})}{\partial r}\right|\Bigg)_{r=1}\\
&=&\lim_{r\rightarrow1-}\frac{|f(e^{i\theta})-f(re^{i\theta})|}{1-r}\\
&\geq&\lim_{r\rightarrow1-}
\frac{1-\frac{1}{2\pi}\int_{0}^{2\pi}K_{\alpha}\big(re^{i(\theta-t)}\big)dt+\frac{(1-r)^{\alpha+1}}{1+r}K_{\alpha}(0)}{1-r}\\
&=&\lim_{r\rightarrow1-}
\frac{1-\frac{1}{2\pi}\int_{0}^{2\pi}K_{\alpha}\big(re^{i\eta}\big)d\eta+\frac{(1-r)^{\alpha+1}}{1+r}K_{\alpha}(0)}{1-r}\\
\end{eqnarray*}
\begin{eqnarray*}
&=&\frac{1}{2\pi}\lim_{r\rightarrow1-}\frac{d}{dr}\int_{0}^{2\pi}K_{\alpha}\big(re^{i\eta}\big)d\eta\\
&=&\lim_{r\rightarrow1-}\frac{d}{dr}\mathcal{M}_{\alpha}(r).
\end{eqnarray*} where $\mathcal{M}_{\alpha}(r)$ is given by {\rm (\ref{eq-11})}.
It follows from the proof of \cite[Theorem 3.1]{O} that
$$\mathcal{M}_{\alpha}(r)=\frac{\left[\Gamma\big(1+\frac{\alpha}{2}\big)\right]^{2}}{\Gamma(1+\alpha)}
F\left(-\frac{\alpha}{2},-\frac{\alpha}{2};1;r^{2}\right)=\frac{\left[\Gamma\big(1+\frac{\alpha}{2}\big)\right]^{2}}{\Gamma(1+\alpha)}
\sum_{n=0}^{\infty}\frac{\left[\big(-\frac{\alpha}{2}\big)_{n}\right]^{2}}{(n!)^{2}}r^{2n},$$
which yields that
$$\frac{d}{dr}\mathcal{M}_{\alpha}(r)=\frac{\alpha^{2}}{2}\frac{\left[\Gamma\big(1+\frac{\alpha}{2}\big)\right]^{2}}{\Gamma(1+\alpha)}
rF\left(1-\frac{\alpha}{2},1-\frac{\alpha}{2};2;r^{2}\right),$$
where $r\in(0,1).$ By (\ref{eq-y}), for $\alpha>0$, we see that

\begin{eqnarray*}
\lim_{r\rightarrow1-}\frac{d}{dr}\mathcal{M}_{\alpha}(r)&=&\frac{\alpha^{2}}{2}\frac{\left[\Gamma\big(1+\frac{\alpha}{2}\big)\right]^{2}}{\Gamma(1+\alpha)}
F\left(1-\frac{\alpha}{2},1-\frac{\alpha}{2};2;1\right)\\
&=&\frac{\alpha^{2}}{2}\frac{\left[\Gamma\big(1+\frac{\alpha}{2}\big)\right]^{2}}{\Gamma(1+\alpha)}
\frac{\Gamma(2)\Gamma(\alpha)}{\left[\Gamma\big(1+\frac{\alpha}{2}\big)\right]^{2}}\\
&=&\frac{\alpha^{2}}{2}\frac{\Gamma(\alpha)}{\Gamma(1+\alpha)}=\frac{\alpha}{2}.
\end{eqnarray*} Therefore, for $\theta\in[0,2\pi]$,
$$\|D_{f}(e^{i\theta})\|\geq \lim_{r\rightarrow1-}\frac{d}{dr}\mathcal{M}_{\alpha}(r)=\frac{\alpha}{2},$$
where $\alpha>0.$ The proof of this theorem is complete. \hfill
$\Box$

\subsection*{Proof of Theorem \ref{thm-1}}  For $r\in[0,1),$ let
$$A_{k}(r,\alpha)=c_{k}F\Big(-\frac{\alpha}{2},k-\frac{\alpha}{2};k+1;r^{2}\Big)$$ and $$B_{k}(r,\alpha)
=c_{-k}F\Big(-\frac{\alpha}{2},k-\frac{\alpha}{2};k+1;r^{2}\Big),$$
where $r=|z|.$ Then

$$A_{k}(r,\alpha)r^{k}=\frac{1}{2\pi}\int_{0}^{2\pi}f(z)e^{-ik\theta}d\theta$$
and 
$$B_{k}(r,\alpha)r^{k}=\frac{1}{2\pi}\int_{0}^{2\pi}f(z)e^{ik\theta}d\theta,$$ 
which imply that

 \be\label{eq-7}
|A_{k}(r,\alpha)|r^{k}=\frac{1}{2\pi}\int_{0}^{2\pi}f(z)e^{-ik\theta}e^{-i\arg
A_{k}(r,\alpha)}d\theta\ee
and
\be\label{eq-8}
|B_{k}(r,\alpha)|r^{k}=\frac{1}{2\pi}\int_{0}^{2\pi}f(z)e^{ik\theta}e^{-i\arg B_{k}(r,\alpha)}d\theta\,,
\ee
where
$A_{k}(r,\alpha)=|A_{k}(r,\alpha)|e^{i\arg A_{k}(r,\alpha)}$,
$B_{k}(r,\alpha)=|B_{k}(r,\alpha)|e^{i\arg B_{k}(r,\alpha)}$ and
$z=re^{i\theta}$. By (\ref{eq-7}),  (\ref{eq-8}) and \cite[Lemma
1]{CR}, we have

\beq\label{eq-9}
&&\left| \big( |A_{k}(r,\alpha)| + |B_{k}(r,\alpha)| \big) r^{k} \right|\\
\nonumber &=&\left| \frac{1}{2\pi} \int_{0}^{2\pi}
                    f(z) \Big[e^{-i \big(k\theta+\arg A_{k}(r,\alpha)\big)} +e^{i \big(k\theta-\arg B_{k}(r,\alpha)\big)} \Big] d\theta \right|\\
\nonumber &\leq& \frac{1}{2\pi}\int_{0}^{2\pi}
                   |f(z)| \Big|e^{-i\big(k\theta+\arg A_{k}(r,\alpha)\big)}+e^{i\big(k\theta-\arg B_{k}(r,\alpha)\big)} \Big| d\theta\\
\nonumber
&\leq&\frac{M}{2\pi}\int_{0}^{2\pi}\Big|e^{-i\big(k\theta+\arg
A_{k}(r,\alpha)\big)}+e^{i\big(k\theta-\arg
B_{k}(r,\alpha)\big)}\Big|d\theta\\ \nonumber
&=&\frac{M}{2\pi}\int_{0}^{2\pi}\Big|1+e^{i\big(2k\theta+\arg
A_{k}(r,\alpha)-\arg B_{k}(r,\alpha)\big)}\Big|d\theta\\ \nonumber
&=&\frac{M}{\pi}\int_{0}^{2\pi}\left|\cos\Bigg(k\theta+\frac{\arg
A_{k}(r,\alpha)-\arg B_{k}(r,\alpha)}{2}\Bigg)\right|d\theta\\
\nonumber &=&\frac{4M}{\pi}.\eeq


 By
letting $r\rightarrow1-$ on (\ref{eq-9}), we obtain
$$|A_{k}(1,\alpha)|+|B_{k}(1,\alpha)|\leq\frac{4M}{\pi}.$$

On the other hand, for $k=0$, we have \beq\label{eq-17}
\frac{1}{2\pi}\lim_{r\rightarrow1-}\int_{0}^{2\pi}|f(re^{i\theta})|^{2}d\theta&=&
\Big|c_{0}F\left(-\frac{\alpha}{2},-\frac{\alpha}{2};1;1\right)\Big|^{2}\\
\nonumber
&&+\sum_{k=1}^{\infty}\Bigg(\bigg|c_{k}F\left(-\frac{\alpha}{2},k-\frac{\alpha}{2};k+1;1\right)\bigg|^{2}\\
\nonumber &&+
\bigg|c_{-k}F\left(-\frac{\alpha}{2},k-\frac{\alpha}{2};k+1;1\right)\bigg|^{2}\Bigg)\\
\nonumber &\leq&M^{2}, \eeq where $r\in[0,1)$. It follows from
(\ref{eq-17}) that
$$\Big|c_{0}F\left(-\frac{\alpha}{2},-\frac{\alpha}{2};1;1\right)\Big|\leq
M.$$ 



If $\alpha=0$, then the  sharpness part follows from \cite[Lemma
1]{CPW4}. The proof of this theorem is complete. \hfill $\Box$

\section{ The Landau type theorem }\label{csw-sec3}

\begin{lem}\label{lem-cm} For $x\in[0,1)$, let
$$\varphi(x)=\frac{\delta}{M(2+\alpha)}-\frac{4Mx}{\pi}\Bigg[\frac{(2-x)}{(1-x)^{2}}
+\frac{2x}{(1-x)(1-x^{2})^{2}}\Bigg],$$ where $\alpha>-2$,
$\delta>0$ and $M>0$ are constant. Then $\varphi$ is strictly
decreasing and there is an unique $x_{0}\in(0,1)$ such that
$\varphi(x_{0})=0.$
\end{lem}
\bpf For $x\in[0,1)$, let
$$f_{1}(x)=\frac{4M}{\pi}\frac{x(2-x)}{(1-x)^{2}}~\mbox{and}~f_{2}(x)=\frac{4M}{\pi}\frac{2x^{2}}{(1-x)(1-x^{2})^{2}}.$$
Since, for $x\in[0,1)$,
$$f_{1}'(x)=\frac{8M}{\pi}\frac{1}{(1-x)^{3}}>0,$$ we see that
$f_{1}$ is continuous and strictly increasing in $[0,1)$. We observe
that $f_{2}$ is also continuous and strictly increasing in $[0,1)$.
Then
$$\varphi(x)=\frac{\delta}{M(2+\alpha)}-f_{1}(x)-f_{2}(x)$$ is continuous and strictly
decreasing in $[0,1)$, which, together with
$$\lim_{x\rightarrow0}\varphi(x)=\frac{\delta}{M(2+\alpha)}>0~\mbox{and}~\lim_{x\rightarrow1-}\varphi(x)=-\infty,$$
imply that there is an unique $x_{0}\in(0,1)$ such that
$\varphi(x_{0})=0.$ \epf

\begin{lem}\label{lem-4}
For $\alpha\in(-1,0]$, let $f\in\mathcal{C}^{2}(\mathbb{D})$ be a
solution to  {\rm (\ref{eq-1})} satisfying $f(0)=|J_{f}(0)|-\beta=0$
and $\sup_{z\in\mathbb{D}}|f(z)|\leq M$, where $M,$ $\beta$ are
positive constants and $J_{f}$ is the Jacobian of $f$. Then $f$ is
univalent in $\mathbb{D}_{\rho_{0}}$, where $\rho_{0}$ satisfies the
following equation
$$\frac{\beta}{M(2+\alpha)}-\frac{4M\rho_{0}}{\pi}\Bigg[\frac{2-\rho_{0}}{(1-\rho_{0})^{2}}
+\frac{2\rho_{0}}{(1-\rho_{0})(1-\rho_{0}^{2})^{2}}\Bigg]=0.$$
Moreover, $f\big(\mathbb{D}_{\rho_{0}}\big)$ contains a univalent
disk $\mathbb{D}_{R_{0}}$ with
$$R_{0}\geq\frac{2\rho_{0}}{3}\left[\frac{\beta}{M(2+\alpha)}-\frac{M\rho_{0}(2-\rho_{0})}{\pi(1-\rho_{0})^{2}}\right].
$$
\end{lem}

\bpf By Theorem \Ref{ThmA}, we can assume that  \beq \nonumber
f(z)&=&\sum_{k=0}^{\infty}c_{k}F\left(-\frac{\alpha}{2},k-\frac{\alpha}{2};k+1;|z|^{2}\right)z^{k}\\
\nonumber &&+
\sum_{k=1}^{\infty}c_{-k}F\left(-\frac{\alpha}{2},k-\frac{\alpha}{2};k+1;|z|^{2}\right)\overline{z}^{k},~z\in\mathbb{D},
\eeq for some sequence $\{c_{k}\}_{k=-\infty}^{\infty}$ of complex
numbers satisfying
$$\lim_{|k|\rightarrow\infty}\sup|c_{k}|^{\frac{1}{|k|}}\leq1.$$



For $\alpha\in(-1,0]$, by \cite[Proposition 1.4]{O}, we observe that
$$F\Big(-\frac{\alpha}{2},k-\frac{\alpha}{2};k+1;r^{2}\Big)=\sum_{n=0}^{\infty}
\frac{\left(-\frac{\alpha}{2}\right)_{n}\left(k-\frac{\alpha}{2}\right)_{n}}{(k+1)_{n}}\frac{r^{2n}}{n!}\geq0$$
is bounded and increasing on $r\in[0,1)$, which imply that

\beq\label{eq-20}&&\big(|c_{k}|+|c_{-k}|\big)F\Big(-\frac{\alpha}{2},k-\frac{\alpha}{2};k+1;r^{2}\Big)\\
\nonumber &\leq&
\big(|c_{k}|+|c_{-k}|\big)F\Big(-\frac{\alpha}{2},k-\frac{\alpha}{2};k+1;1\Big)\\
\nonumber &\leq&\frac{4M}{\pi}, \eeq where $r=|z|$ and
$k\in\{1,2,\ldots\}$.

By (\ref{eq-20}) and Theorem \ref{thm-1}, we see that, for each
$k\in\{1,2,\ldots\}$, \be\label{eq-16}
(|c_{k}|+|c_{-k}|)\frac{\left(-\frac{\alpha}{2}\right)_{n}
\left(k-\frac{\alpha}{2}\right)_{n}}{(k+1)_{n}}\frac{1}{n!}\leq\frac{4M}{\pi},\ee
where $n\in\{1,2,\ldots\}$.

Since  $c_{0}=f(0)=0$, we see that

\beq\label{eq-14}
f_{z}(z)-f_{z}(0)&=&\sum_{k=2}^{\infty}kc_{k}F\left(-\frac{\alpha}{2},k-\frac{\alpha}{2};k+1;w\right)z^{k-1}\\
\nonumber
&&+\sum_{k=1}^{\infty}c_{k}\frac{d}{dw}F\left(-\frac{\alpha}{2},k-\frac{\alpha}{2};k+1;w\right)z^{k}\overline{z}\\
\nonumber
&&+\sum_{k=1}^{\infty}c_{-k}\frac{d}{dw}F\left(-\frac{\alpha}{2},k-\frac{\alpha}{2};k+1;w\right)\overline{z}^{k+1}
\eeq and \beq\label{eq-15}
f_{\overline{z}}(z)-f_{\overline{z}}(0)&=&\sum_{k=2}^{\infty}kc_{-k}F\left(-\frac{\alpha}{2},k-\frac{\alpha}{2};k+1;w\right)\overline{z}^{k-1}\\
\nonumber
&&+\sum_{k=1}^{\infty}c_{k}\frac{d}{dw}F\left(-\frac{\alpha}{2},k-\frac{\alpha}{2};k+1;w\right)z^{k+1}\\
\nonumber
&&+\sum_{k=1}^{\infty}c_{-k}\frac{d}{dw}F\left(-\frac{\alpha}{2},k-\frac{\alpha}{2};k+1;w\right)\overline{z}^{k}z,
\eeq where $w=|z|^{2}.$


Applying   (\ref{eq-16}), (\ref{eq-14}) and (\ref{eq-15}), we obtain

\beq\label{eq-18}
&&|f_{z}(z)-f_{z}(0)|+|f_{\overline{z}}(z)-f_{\overline{z}}(0)|\\
\nonumber &\leq&
\sum_{k=2}^{\infty}k\big(|c_{k}|+|c_{-k}|\big)F\left(-\frac{\alpha}{2},k-\frac{\alpha}{2};k+1;w\right)|z|^{k-1}\\
\nonumber
&&+2\sum_{k=1}^{\infty}\big(|c_{k}|+|c_{-k}|\big)\frac{d}{dw}F\left(-\frac{\alpha}{2},k-\frac{\alpha}{2};k+1;w\right)|z|^{k+1}\\
 \nonumber
&\leq&\frac{4M}{\pi}\sum_{k=2}^{\infty}k|z|^{k-1}+2\sum_{k=1}^{\infty}\left[\frac{4M}{\pi}\sum_{n=1}^{\infty}n|z|^{2(n-1)}\right]|z|^{k+1}\\
\nonumber&=&\frac{4M}{\pi}\frac{|z|(2-|z|)}{(1-|z|)^{2}}+\frac{8M}{\pi}\sum_{k=1}^{\infty}\frac{|z|^{k+1}}{(1-|z|^{2})^{2}}
\\ \nonumber
&=&\frac{4M}{\pi}\frac{|z|(2-|z|)}{(1-|z|)^{2}}+\frac{8M}{\pi}\frac{|z|^{2}}{(1-|z|)(1-|z|^{2})^{2}}.
\eeq

Applying Theorem \ref{thm-2} (\ref{eq-0.12}), we get
$$\beta=|J_{f}(0)|=|\det
D_{f}(0)|=\|D_{f}(0)\|l\big(D_{f}(0)\big)\leq
M(2+\alpha)l\big(D_{f}(0)\big),$$ which gives that
\be\label{eq-0.14}l\big(D_{f}(0)\big)\geq\frac{\beta}{M(2+\alpha)}.
\ee

In order to prove the univalence of $f$ in $\mathbb{D}_{\rho_{0}}$,
we choose two distinct points $z_{1}$, $z_{2}\in
\mathbb{D}_{\rho_{0}}$ and let $[z_{1}, z_{2}]$ denote the segment
from $z_{1}$ to $z_{2}$ with the endpoints $z_{1}$  and $z_{2}$,
where $\rho_{0}$ satisfies the following equation
$$\frac{\beta}{M(2+\alpha)}-\frac{4M\rho_{0}}{\pi}\Bigg[\frac{2-\rho_{0}}{(1-\rho_{0})^{2}}
+\frac{2\rho_{0}}{(1-\rho_{0})(1-\rho_{0}^{2})^{2}}\Bigg]=0.$$

By (\ref{eq-18}), (\ref{eq-0.14}) and Lemma \ref{lem-cm}, we have

\begin{eqnarray*}
|f(z_{2})-f(z_{1})|&=&\left|\int_{[z_{1},z_{2}]}f_{z}(z)dz+f_{\overline{z}}(z)d\overline{z}\right|\\
&=&\left|\int_{[z_{1},z_{2}]}f_{z}(0)dz+f_{\overline{z}}(0)d\overline{z}\right|\\
&&-\left|\int_{[z_{1},z_{2}]}\big(f_{z}(z)-f_{z}(0)\big)dz+\big(f_{\overline{z}}(z)-f_{\overline{z}}(0)\big)d\overline{z}\right|\\
&\geq&l(D_{f})(0)|z_{2}-z_{1}|\\
&&-\int_{[z_{1},z_{2}]}\big(|f_{z}(z)-f_{z}(0)|+|f_{\overline{z}}(z)-f_{\overline{z}}(0)|\big)|dz|\\
\end{eqnarray*}
\begin{eqnarray*}
&>&|z_{2}-z_{1}|\Bigg\{\frac{\beta}{M(2+\alpha)}\\
&&-\frac{4M\rho_{0}}{\pi}\Bigg[\frac{2-\rho_{0}}{(1-\rho_{0})^{2}}
+\frac{2\rho_{0}}{(1-\rho_{0})(1-\rho_{0}^{2})^{2}}\Bigg]\Bigg\}\\
&=&0.
\end{eqnarray*}
Thus, $f(z_{2})\neq f(z_{1}).$ The univalence of $f$ follows from
the arbitrariness of $z_{1}$ and $z_{2}$. This implies that $f$ is
univalent in $\mathbb{D}_{\rho_{0}}$.

Now, for any $\zeta'=\rho_{0}
e^{i\theta}\in\partial\mathbb{D}_{\rho_{0}}$, we  obtain that
\begin{eqnarray*}
|f(\zeta')-f(0)|&=&\left|\int_{[0,\zeta']}f_{z}(z)dz+f_{\overline{z}}(z)d\overline{z}\right|\\
&=&\left|\int_{[0,\zeta']}f_{z}(0)dz+f_{\overline{z}}(0)d\overline{z}\right|\\
&&-\left|\int_{[0,\zeta']}\big(f_{z}(z)-f_{z}(0)\big)dz+\big(f_{\overline{z}}(z)-f_{\overline{z}}(0)\big)d\overline{z}\right|\\
&\geq&l(D_{f})(0)\rho_{0}-\int_{[0,\zeta']}\big(|f_{z}(z)-f_{z}(0)|+|f_{\overline{z}}(z)-f_{\overline{z}}(0)|\big)|dz|\\
&\geq&l(D_{f})(0)\rho_{0}-\frac{4M\rho_{0}^{2}}{\pi}\int_{0}^{1}\left[\frac{t(2-\rho_{0}t)}{(1-\rho_{0}t)^{2}}+
\frac{2\rho_{0}t^{2}}{(1-\rho_{0}t)(1-\rho_{0}^{2}t^{2})^{2}}\right]dt\\
&\geq&\frac{\beta\rho_{0}}{M(2+\alpha)}-\frac{4M\rho_{0}^{2}}{\pi}\Bigg[\frac{(2-\rho_{0})}{(1-\rho_{0})^{2}}\int_{0}^{1}t
dt\\
&&+\frac{2\rho_{0}}{(1-\rho_{0})(1-\rho_{0}^{2})^{2}}\int_{0}^{1}t^{2}dt\Bigg]\\
&=&\rho_{0}\left\{\frac{\beta}{M(2+\alpha)}-\frac{4M\rho_{0}}{\pi}\left[\frac{2-\rho_{0}}{2(1-\rho_{0})^{2}}+
\frac{2\rho_{0}}{3(1-\rho_{0})(1-\rho_{0}^{2})^{2}}\right]\right\}\\
&=&\frac{2\rho_{0}}{3}\left[\frac{\beta}{M(2+\alpha)}-\frac{M\rho_{0}(2-\rho_{0})}{\pi(1-\rho_{0})^{2}}\right].
\end{eqnarray*}
Hence $f\big(\mathbb{D}_{\rho_{0}}\big)$ contains a univalent disk
$\mathbb{D}_{R_{0}}$ with
$$R_{0}\geq\frac{2\rho_{0}}{3}\left[\frac{\beta}{M(2+\alpha)}-\frac{M\rho_{0}(2-\rho_{0})}{\pi(1-\rho_{0})^{2}}\right]. $$
The proof of this lemma is complete. 
\end{pf}

Let us recall the following result which is referred to as {\it
Jensen's inequality} (cf. \cite{R}).

\begin{Lem}\label{Jensen}
Let $(\Omega, A, \mu)$ be a measure space such that $\mu(\Omega)=1.$
If $g$ is a real-valued function that is $\mu$-integrable, and if
$\chi$ is a convex function on the real line, then
$$\chi\left (\int_{\Omega}g\,d\mu\right )\leq\int_{\Omega}\chi\circ g\,d\mu.
$$
\end{Lem}

\subsection*{Proof of Theorem \ref{thm-4}}
For $z\in\mathbb{D}_{r}$, we have $$f(z)=\frac{c_{\alpha}}{2\pi
r^{\alpha}}\int_{0}^{2\pi}\frac{\big(r^{2}-|z|^{2}\big)^{\alpha+1}}{|r-ze^{-it}|^{2+\alpha}}f(re^{it})dt,$$
where $r\in(0,1)$. Let
$$\phi_{z}(r)=\frac{c_{\alpha}}{2\pi
r^{\alpha}}\int_{0}^{2\pi}\frac{\big(r^{2}-|z|^{2}\big)^{\alpha+1}}{|r-ze^{-it}|^{2+\alpha}}dt,$$
where $z\in\mathbb{D}_{r}$. Applying \cite[Theorem 3.1]{O}, we see
that, for $z\in\mathbb{D}$,
\be\label{eq-22}\phi_{z}(1)\leq\lim_{|z|\rightarrow1-}\phi_{z}(1)=1.\ee
By using Jensen's inequality (see Lemma \Ref{Jensen}), for $p\geq1$,
we get

\begin{eqnarray*}
\left|\frac{f(z)}{\phi_{z}(r)}\right|^{p}&=&\left|\frac{1}{2\pi}\int_{0}^{2\pi}\frac{c_{\alpha}}{r^{\alpha}\phi_{z}(r)}
\frac{\big(r^{2}-|z|^{2}\big)^{\alpha+1}}{|r-ze^{-it}|^{2+\alpha}}f(re^{it})dt\right|^{p}\\
&\leq&\frac{1}{2\pi}\int_{0}^{2\pi}\frac{c_{\alpha}}{r^{\alpha}\phi_{z}(r)}
\frac{\big(r^{2}-|z|^{2}\big)^{\alpha+1}}{|r-ze^{-it}|^{2+\alpha}}|f(re^{it})|^{p}dt\\
&\leq&\frac{c_{\alpha}}{r^{\alpha}\phi_{z}(r)}\frac{\big(r^{2}-|z|^{2}\big)^{\alpha+1}}{\big(r-|z|\big)^{2+\alpha}}
\left(\frac{1}{2\pi}\int_{0}^{2\pi}|f(re^{it})|^{p}dt\right)\\
&\leq&\frac{c_{\alpha}\|f\|_{p}^{p}}{r^{\alpha}\phi_{z}(r)}\frac{\big(r+|z|\big)^{\alpha+1}}{\big(r-|z|\big)},
\end{eqnarray*}
which implies that
$$|f(z)|\leq\left[\frac{c_{\alpha}\|f\|_{p}^{p}\big(\phi_{z}(r)\big)^{p-1}}{r^{\alpha}}\right]^{\frac{1}{p}}
\frac{\big(r+|z|\big)^{\frac{\alpha+1}{p}}}{\big(r-|z|\big)^{\frac{1}{p}}},$$
where  $z\in\mathbb{D}_{r}$. By letting $r\rightarrow1-$ and
(\ref{eq-22}), for $z\in\mathbb{D},$ we have
\be\label{eq-23}|f(z)|\leq\left[c_{\alpha}\|f\|_{p}^{p}\big(\phi_{z}(1)\big)^{p-1}\right]^{\frac{1}{p}}
\frac{\big(1+|z|\big)^{\frac{\alpha+1}{p}}}{\big(1-|z|\big)^{\frac{1}{p}}}\leq
c_{\alpha}^{\frac{1}{p}}\|f\|_{p}\frac{\big(1+|z|\big)^{\frac{\alpha+1}{p}}}{\big(1-|z|\big)^{\frac{1}{p}}}.\ee
For $\zeta\in\mathbb{D},$  let $Q(\zeta)=f(\gamma\zeta)/\gamma,$
where $\gamma\in(0,1)$. It is not difficult to know that
$Q(0)=|J_{Q}(0)|-\lambda=0$.  By (\ref{eq-23}), for
$\zeta\in\mathbb{D},$ we obtain

$$|Q(\zeta)|=\frac{|f(\gamma\zeta)|}{\gamma}\leq c_{\alpha}^{\frac{1}{p}}\|f\|_{p}
\frac{\big(1+\gamma\big)^{\frac{\alpha+1}{p}}}{\gamma\big(1-\gamma\big)^{\frac{1}{p}}},$$
which gives that $$|Q(\zeta)|\leq
c_{\alpha}^{\frac{1}{p}}\|f\|_{p}\min_{0<\gamma<1}\mu(\gamma),$$
where
$$\mu(\gamma)=\big(1+\gamma\big)^{\frac{\alpha+1}{p}}/
\left[\gamma\big(1-\gamma\big)^{\frac{1}{p}}\right].$$ Let
$\gamma_{0}\in(0,1)$ satisfy
$$\mu(\gamma_{0})=\min_{0<\gamma<1}\mu(\gamma).$$ By using Lemma
\ref{lem-4}, we observe that $Q$ is univalent in
$\mathbb{D}_{\rho_{0}}$, where $\rho_{0}$ satisfies the following
equation
$$\frac{\lambda}{M^{\ast}(2+\alpha)}-\frac{4M^{\ast}\rho_{0}}{\pi}\Bigg[\frac{2-\rho_{0}}{(1-\rho_{0})^{2}}
+\frac{2\rho_{0}}{(1-\rho_{0})(1-\rho_{0}^{2})^{2}}\Bigg]=0,$$ where
$M^{\ast}=c_{\alpha}^{\frac{1}{p}}\|f\|_{p}\mu(\gamma_{0}).$
Moreover, $Q\big(\mathbb{D}_{\rho_{0}}\big)$ contains a univalent
disk $\mathbb{D}_{R_{0}}$ with
$$R_{0}\geq\frac{2\rho_{0}}{3}\left[\frac{\lambda}{M^{\ast}(2+\alpha)}-\frac{M^{\ast}\rho_{0}(2-\rho_{0})}{\pi(1-\rho_{0})^{2}}\right].
$$ Hence $f$ is univalent in
$\mathbb{D}_{\gamma_{0}\rho_{0}}$ and
$f\big(\mathbb{D}_{\gamma_{0}\rho_{0}}\big)$ contains a univalent
disk $\mathbb{D}_{\gamma_{0}R_{0}}$. The proof of this theorem is
complete. \hfill $\Box$

\bigskip

{\bf Acknowledgements:} This research was partly  supported by the
National Natural Science Foundation of China (No. 11401184 and No.
11326081), the Hunan Province Natural Science Foundation of China
(No. 2015JJ3025), the Excellent Doctoral Dissertation of Special
Foundation of Hunan Province (higher education 2050205), the
Construct Program of the Key Discipline in Hunan Province, the
V\"ais\"al\"a Foundation of The Finnish Academy of Sciences and
Letters.





\normalsize


\begin{thebibliography}{99}

\bibitem{AA} {\sc Z. Abdulhadi and Y. Abu Muhanna,}
Landau's theorem for biharmonic mappings, \textit{J. Math. Anal.
Appl.,} {\bf 338}(2008), 705--709.

\bibitem{AP} {\sc A. Alemam and J. A. Pel\'aez,}
Spectra of integration operators and weighted square functions,
\textit{Indiana Univ. Math. J.,} {\bf 61}(2012), 775--793.

\bibitem{ABRVV} {\sc G. D. Anderson, R. W. Barnard, K. C. Richards, M. K. Vamanamurthy and M. Vuorinen,}
Inequalities for zero-balanced hypergeometric functions,
\textit{Trans. Amer. Math. Soc.}, {\bf 126}(1995), 1713--1723.

\bibitem{QVV} {\sc G. D. Anderson, S. L. Qiu, M. K. Vamanamurthy and M. Vuorinen,} Generalized elliptic integrals and modular equations,
\textit{Pacific J. Math.}, {\bf 192}(2000), 1--37.

\bibitem{AH} {\sc A. Borichev and H. Hedenmalm,}
Weighted integrability of polyharmonic functions, \textit{Adv.
Math.,} {\bf 264}(2014), 464--505.


\bibitem{H1} {\sc H. H. Chen,} The Schwarz-Pick lemma for planar harmonic
mappings, \textit{ Sci. China Math.}, {\bf 54}(2011), 1101--1118.


\bibitem{HG} {\sc H. Chen, P. M. Gauthier and W. Hengartner,}
Bloch constants for planar harmonic mappings, \textit{Proc. Amer.
Math. Soc.,} {\bf 128}(2000), 3231--3240.

\bibitem{CPW1} {\sc Sh. Chen, S.~Ponnusamy and A. Rasila,} Coefficient estimates, Landau¡¯s theorem and Lipschitz-type spaces on planar harmonic
mappings, \textit{J. Aust. Math. Soc.,} {\bf 96}(2014), 198--215.

\bibitem{CPW2} {\sc Sh. Chen, S. Ponnusamy and X. Wang,} Harmonic mappings in Bergman spaces,
\textit{Monatsh. Math.,} {\bf 170}(2013), 325--342.

\bibitem{CPW3} {\sc Sh. Chen, S. Ponnusamy and X. Wang,} Integral means and coefficient estimates on planar harmonic
mappings, \textit{Ann. Acad. Sci. Fenn. Math.,} {\bf 37}(2012),
69--79.

\bibitem{CPW5} {\sc Sh. Chen, S. Ponnusamy and X. Wang,} On planar harmonic Lipschitz and planar harmonic Hardy classes,
 \textit{Ann. Acad. Sci. Fenn. Math.,} {\bf 36}(2011),
567--576.



\bibitem{CPW4} {\sc Sh. Chen, S. Ponnusamy and X. Wang,}
Bloch constant and Landau's theorems for planar \textit{p-}harmonic
mappings, \textit{J. Math. Anal. Appl.,} {\bf 373}(2011), 102--110.


\bibitem{CR} {\sc Sh. Chen and A. Rasila,} Schwarz-Pick type estimates of pluriharmonic mappings in the unit
polydisk,   arXiv:1409.7897, to appear in \textit{Illinois J. Math.}

\bibitem{Co} {\sc F.~ Colonna,}
The Bloch constant of bounded harmonic mappings, \textit{Indiana
Univ. Math. J.,} {\bf 38} (1989), 829--840.

\bibitem{Du} {\sc P. Duren,}
\textit{Harmonic Mappings in the Plane,} Cambridge Univ. Press,
2004.

\bibitem{E}  {\sc M. Engli$\breve{{\rm s}}$,}
A Loewner-type lemma for weighted biharmonic operators,
\textit{Pacific J. Math.}, {\bf 179}(1997), 343--353.

\bibitem{HZ}  {\sc E. Heinz,}
On one-to-one harmonic mappings, \textit{Pacific J. Math.}, {\bf
9}(1959), 101--105.



\bibitem{K}  {\sc D. Kalaj,}
Heinz inequality for the unit ball,  arXiv: 1504.01686 [math. AP].


\bibitem{KV}  {\sc D. Kalaj and M. Vuorinen,}
On harmonic functions and the Schwarz lemma, \textit{Proc. Amer.
Math. Soc.}, {\bf 140}(2012), 161--165.

\bibitem{L} {\sc E. Landau,} \"Uber die Bloch'sche konstante und zwei verwandte weltkonstanten,
\textit{Math. Z.,} {\bf 30}(1929), 608--634.

\bibitem{LZ} {\sc M. S. Liu, Z. W. Liu and Y. C. Zhu,}
Landau's theorems for certain biharmonic mappings, \textit{Acta
Math. Sinica, Chinese Series,} {\bf 54}(2011), 1--12.


\bibitem{O1} {\sc A. Olofsson,}
A representation formula for radially weighted biharmonic functions
in the unit disc, \textit{Publ. Math.,} {\bf 49} (2005), 393--415.

\bibitem{O2} {\sc A. Olofsson,}
Regularity in a singular biharmonic Dirichlet problem,
\textit{Monatsh. Math.,} {\bf 148} (2006), 229--239.




\bibitem{O} {\sc A. Olofsson,}
Differential operators for a scale of Poisson type kernels in the
unit disc, \textit{J. Anal. Math.,} {\bf 123} (2014), 227--249.


\bibitem{O3} {\sc A. Olofsson and J. Wittsten,}
 Poisson integrals for standard weighted
Laplacians  in the unit disc, \textit{J. Math. Soc. Japan.,} {\bf
65} (2013), 447--486.



 \bibitem{P1} {\sc M. Pavlovi\'c,} Decompositions of $L^{p}$   and Hardy spaces of polyharmonic
 functions,
 \textit{ J. Math. Anal. Appl.,}  {\bf 216} (1997), 499--509.


\bibitem{P0} {\sc M. Pavlovi\'c,} Harmonic Schwarz lemmas: Chen.
Kalaj-Vuorinen. Pavlovi\'c. Heinz., Preprint.











 \bibitem{R} {\sc W. Rudin,}  Real and Complex Analysis. McGraw-Hill. ISBN
0-07-054234-1,1987.

\bibitem{SW} {\sc W. Szapiel,} Bounded harmonic mappings,
\textit{J. Anal. Math.,} {\bf 111}(2010), 47--76.

\end{thebibliography}
\end{document}